\def\proof{{\noindent{\bf Proof:} \quad }}
\def\endproof{\hfill{\vrule height4pt width6pt depth2pt} \vskip .2in}

\def\to{\rightarrow}

\def\vec3#1#2#3{(#1,#2,#3)}
\def\emp{\emptyset}

\def\ced{\c}

\def\proof{{\noindent{\bf Proof:} \qquad }}

\def\set#1{\{#1\}}

\def\al{\alpha}

\def\de{\delta}
\def\ep{\epsilon}

\def\th{\theta}

\def\si{\sigma}

\def\ph{\varphi}

\def\ps{\psi}

\def\today{{\ifcase\month\or January\or February\or March \or April
\or May\or June\or July \or August \or September \or
October \or November \or December\fi}\  \number\day,\
\number\year}


\def\su #1 #2 #3 {#1_{#2}^{#3}}

\def\abs#1{|#1|}

\def\X{(X, {\cal B}, \mu,T)}
\def\Y{(Y, {\cal C}, \nu, S)}
\def\N{\Bbb N}
\def\R{\Bbb R}
\def\Z{\Bbb Z}
\def\C{\Bbb C}

\def\T{\Bbb T}

\def\emptyset{\varnothing}




\magnification=1300     

\baselineskip=16pt           
   \input amssym.tex

   \def\C{{\cal C}}               
   \def\xc{(X,\mu,T)}             
   \def\yc{(Y,\nu,S)}            
\hfuzz.1in

\centerline{ Dye's theorem in the almost continuous category}
\vskip.5in
\centerline{ Andr\'es del Junco{$^{1}$}
   and Ay\ced se \ced Sahin$^{2}$}

\vskip.5in

\centerline{ Abstract}
\bigskip

Suppose $X$ and $Y$ are Polish spaces with 
non-atomic Borel probability
measures $\mu$ and $\nu$ and suppose that $T$ and $S$ are ergodic 
measure-preserving homeomorphisms of $(X,\mu)$ and $(Y,\nu)$.
Then there are invariant $G_{\de}$ subsets $X' \subset X$ and $Y' 
\subset Y$ of full measure and a homeomorphism $\ph: X' \to Y'$ which 
maps $\mu|_{X'}$ to $\nu|_{Y'}$ and maps $T$-orbits onto $S$-orbits.
We also deal with the case where $T$ and $S$ preserve infinite 
invariant measures.

\vskip1.5in
\item{$^{1}$} Dept. of Mathematics,  University of Toronto\hfill \break
	     {\sl deljunco@math.toronto.edu}
\medskip
\item{$^{2}$} Dept. of Mathematics, DePaul University, \hfill \break
	      {\sl asahin@condor.depaul.edu}

\vfill\eject

\noindent{Section 1. Introduction}

\bigskip
Suppose $\X$ and $\Y$ are measure-preserving dynamical systems on
probability spaces.  An
orbit equivalence between $X$ and $Y$ is a measurable map
$\ph:X' \to Y'$, where $X'$ and $Y'$ are invariant subsets of measure
one, with a measurable inverse $\ps: Y' \to X'$ such that
$\ph^{*}\mu = \nu$ and $\ph$ maps
$T$-orbits  onto $S$-orbits.  A fundamental theorem of Dye
[D1] asserts that any two measure-preserving systems are orbit
equivalent.

Recently Hamachi and Keane [HK] proved that the binary and ternary
odometers are actually orbit equivalent in a stronger sense, namely
the orbit equivalence $\ph$ can be chosen to be continuous with a
continuous inverse when 
restricted to an appropriate invariant subset of measure $1$.
Such an orbit equivalence $\ph$ is called almost continuous or 
finitary, the term finitary being more appropriate
 in the context of maps defined on sequence spaces,
such as the odometers. We will use the term almost continuous as we
will be working in a more general setting.
There have been a number of results since then asserting the 
existence of an almost continuous orbit equivalence for various 
classes of ergodic measure-preserving homeomorphisms or maps
which are homeomorphisms after restriction to an invariant subset of
full measure, as is the case with odometers.
See [HKR], [R1], [R2] and [RR3].
Perhaps the most general of these is the
result of Hamachi, Keane and Roychowdhury [HKR] asserting that any two
adic transformations are almost continuously orbit equivalent.
We mention here also the celebrated work of Keane and Smorodinsky 
[KS1], [KS2] on finitary isomorphism, the paper of 
del Junco [J] on finitary unilateral isomorphism and the recent papers of 
Roychowdhury and Rudolph [RR1], [RR2] on finitary Kakutani equivalence.

\def\ot{{\cal O}_{T}}

Here we will prove the following  general result, an almost continuous 
version of Dye's theorem, which includes all the orbit equivalence results
mentioned above and  much more. 
Recall that a subset of a Polish space is Polish if
and only if it is a $G_\de$.

\proclaim Theorem 1. 
Suppose $X$ and $Y$ are Polish spaces with 
non-atomic Borel probability
measures $\mu$ and $\nu$ and suppose that $T$ and $S$ are ergodic 
measure-preserving homeomorphisms of $(X,\mu)$ and $(Y,\nu)$.
Then there are invariant $G_{\de}$ subsets $X' \subset X$ and $Y' 
\subset Y$ of full measure and a homeomorphism $\ph: X' \to Y'$ which 
maps $\mu|_{X'}$ to $\nu|_{Y'}$ and maps $T$-orbits onto $S$-orbits.

We will also prove the analogous result in the case where the
invariant measure is infinite, Theorem 3 in Section 3 below.
We thank Sasha Danilenko for an
insight which simplified the proof of Theorem 3.

Theorem 1 says that the restrictions of $T$ and $S$ to $X'$ and $Y'$ 
are topologically orbit equivalent, in the sense of Giordano, Putnam 
and Skau [GPS], via a map which also carries $\mu$ to $\nu$.
Our proof is a combination of the techniques used in   well-known proofs 
of Dye's theorem (see for example  [KW] or [HIK])
with the Hamachi-Keane technique.
For those familiar with [HK] we remark that
 we  use constructs very similar to theirs but in a different
 setting with different terminology.

Of course, Theorem 1 also applies to discontinuous $T$ and $S$ as well,
provided they have restrictions to invariant 
dense $G_{\de}$ subsets of full measure which are continuous. This
class includes odometers, interval exchange maps and 
 adic transformations among others.

We remark that if $\mu$ and $\nu$ have full supports then
$X'$ and $Y'$ in Theorem 1
are necessarily dense $G_{\de}$'s.
Theorem 1 is both measure-theoretic and topological in nature.
By the above remark the purely topological character is
that of generic orbit
 equivalence as defined by
Sullivan, Weiss and 
Wright, [SWWr]. They prove a very general result, namely
any two  discrete groups of 
 homeomorphisms of  Polish spaces are continuously orbit equivalent
 after restriction to invariant dense $G_{\de}$ subsets.
 In the special case of single homeomorphisms which possess some
 invariant probability of full support their result follows from
 ours.

Before we proceed with the proof of Theorem 1 we mention some
 questions which arise from this work. First, it is likely that 
our methods will show that any ergodic  action of a discrete
amenable group $G$ by
 homeomorphisms of a Polish space 
 preserving a probability measure is orbit equivalent to such an 
 action with group $G=Z$. 
 Is there an almost continuous analogue to the
 theorem of Dye [D2] which states that any isomorphism of the full
 groups of two countable discrete groups of
 transformations preserving a probability measure  is 
 implemented by an orbit equivalence?
 Is there a theory in the case of non-singular measures with some,
 or all 
 the features of the measurable theory (see [K], [HO], [KW])? 
Is there anything that distinguishes the almost continuous 
classification from the measurable classification?
In other works
might it be the case that if
 two non-singular homeomorphisms 
(or more generally groups of homeomorphisms)
  are measurably orbit equivalent then they
 must also be almost continuously orbit equivalent?
 
 Returning to integer actions, note that
 every orbit equivalence between $T$ and $S$
 is also an isomorphism between $T$ and $S'$ where 
 $S'$ is a map with the same orbits as $S$ so that  $S'y = S^{n(y)}y$
 and $S(y) = S'^{m(y)}y$ where
 $n$  and $m$ are integer valued functions (co-cycles) on $Y$.
 Can we require
 in Theorem 1
 that $m$ and $n$ be continuous on the set $Y'$? This question 
 is motivated by
 the definition due to Giordano, Putnam and Skau [GPS] of
 strong topological orbit equivalence.
 
\bigskip
\noindent{Section 2. Proof of Theorem 1}
 
 To prove Theorem 1 we first reduce  it to a
special case. Let us say a topological space $X$ is
{\bf fractured}
if its topology
has a countable base consisting of clopen subsets.
(We do not know if any zero-dimensional Polish space is
fractured.)
In a fractured space every
 open set is a union of countably many disjoint
clopen sets.
Note that any subset of a 
fractured space is  trivially fractured.

\proclaim Theorem 2. Theorem 1 holds with the additional hypothesis 
that $X$ and $Y$ are fractured.

 How does Theorem 1 follow from Theorem 2? It suffices to 
 show that $X$ has an invariant
  $G_{\de}$ subset $X'$ of measure 1 which is fractured, since
  $X'$ is again Polish.
 To construct $X'$ 
 first observe that for a fixed $x \in X$ at most countably many 
 of the metric spheres $S(x,r)=\set{y: d(x,y) = r}$ can have
positive 
 measure. This means that we may find a countable dense subset
 $ R_x \subset \R^{+}$ such that the spheres $S(x,r), r \in  R_{x}$, all 
 have measure zero.
Let $ \set{x_i: i \in I}$ be a countable dense subset of $X$ and set
$$
X^{\sharp} = (\bigcup_{i \in I}\bigcup_{r \in R_{x_i}}S(x_{i},r))^c .
$$
 $X^\sharp$ is a $G_{\de}$ of measure one and it is easy to see that
$X^{\sharp}$ is fractured, so $X' = \bigcap_{n \in \Z}T^n X^\sharp$ is 
the desired invariant $G_\de$ of measure one.

 \bigskip
 
  \def\T{{\cal T}}               
   \def\S{{\cal S}}  
 
\def\nec{non-empty clopen set }
\def\necs{non-empty clopen sets }

The rest of this section is devoted to the proof of Theorem 2.
We will implicitly use  fixed complete metrics on $X$ and $Y$.
We will state definitions and lemmas for $\xc$ with the
understanding that we will use them for $\yc$ as well.
 Since the support of $\mu$ is an invariant $G_{\de}$ subset of full measure
 we can and shall assume henceforth that the measures  $\mu$ and
 $\nu$ have full support. Thus any non-empty open set has non-zero measure.
The fact that $X$ is fractured and $\mu$ is non-atomic implies 
that every non-empty open set
contains \necs of measure as small as we please.
We will frequently use
(implicitly) the fact that the clopen sets form an algebra invariant 
under $T$ so that all finite set operations involving clopen sets and
powers of $T$ yield clopen sets. If $A$ is a measurable subset of $X$
and $A_{1}, A_{2}, \ldots $ are subsets of $A$ we will say that
$\set{A_{i}}$ {\bf fills out} $A$ if 
$\mu(\bigcup_{i}A_{i}) = \mu(A)$.

\proclaim Lemma 1. Suppose $A$ and $B$ are open subsets of $X$ such that 
$\mu(A) = \mu(B)$. Then there are disjoint 
clopen sets $A_1,A_2, \ldots \subset A$ which fill out $A$ and 
integers $n_1, n_2, \ldots$ such that the sets $B_{i} := T^{n_{i}}A_{i}$
are disjoint and 
contained in $B$ (and fill out $B$).

\proof
Suppose $A_{i}$ and $B_{i}$ have been defined for all $i < m$. Set
$A^{m}= A\backslash \bigcup_{i < m} A_{i}$ and $B^{m}= B\backslash 
\bigcup_{i < m} B_{i}$. Since $A^{m}$ and $B^{m}$ are open we may find
clopen subsets
$A_{m}' \subset A^{m}$ and $B_{m}' \subset B^{m}$ such that  $\mu(A_{m}') > 
\mu(A^{m})/2$ and  $\mu(B_{m}') > 
\mu(B^{m})/2$. By ergodicity of $T$ we may then choose $n_{m}$ such 
that 
$$ \mu(T^{n_{m}}A_{m}' \cap B_{m}') \geq
{1 \over 2}\mu(A_{m}')\mu(B_{m}')
\geq {1 \over 8}\mu(A^{m})\mu(B^{m}).   \eqno (1)$$
We set $B_{m} =T^{n_{m}}A_{m}' \cap B_{m}' $ and $A_{m}= 
T^{-n_{m}}B_{m}$.

To see that the $A_{i}$ fill out A, we must show  that 
 $ \alpha := \lim \mu(A^{n}) = 0$. Since $\mu(A^{n}) = \mu(B^{n}) \geq 
\al$ we will have $\mu(A_{n+1}) \geq \al^{2}/8 $ for all $n$. 
It follows that $\al = 0$ since the $A_{n}$ are disjoint.

\proclaim Lemma 2.
\smallskip
\item{(a)} If $A$ is a  clopen set in $X$ and $\ep > 0$ then $A$ can be 
partitioned into clopen subsets $A_{1}, \ldots, A_{k}$ such that
$\mu(A_{i}) < \ep$. Consequently  any open set can
be written as a countable 
disjoint union of clopen sets with measure less than $\ep$.
\smallskip
\item{(b)} If $B$ is an open set in $X$ and
$\mu(B)= r_{1}+ r_{2 }+ \ldots $ then there are disjoint open subsets 
$B_{1}, B_{2}, \ldots $ of $B$ such that
$\mu(B_{i}) = r_{i}$.

\proof For (a)  let $C$ be any \nec of measure less than $\ep$.
Suppose $A_{1},A_{2}, \ldots, A_{i}$ have been defined. Let
$A^{i} = A \backslash \bigcup_{j=1}^{i}A_{i}$. If $\mu(A^{i}) < \ep$ we 
are done. Otherwise find an $n$ such that 
$\mu(T^{n}C \cap A^{i}) > \mu(C)\mu(A^{i})/2$ and set 
$A_{i+1} = T^{n}C \cap A^{i}$.
 Since we always remove a 
fraction at least $\mu(C)/2$ from  $A^{i}$ we will eventually arrive 
at an $A^{k}$ of measure less than $\ep$ at which point we are done.

\def\emp{\emptyset}

To prove (b) we will construct for each
$n$ a sequence $B_{1}^{n}, B_{2}^{n}, \ldots$ 
of disjoint clopen sets such that $B_{i}^{n} = \emp$ for all 
sufficiently large $i$ and
such that for each $i$ the sequence $B^{1}_{i}, B^{2}_{i}, \ldots$ increases to
the desired $B_{i}$.
Suppose that for a given $n$
we have already found such disjoint clopen subsets
$B_{1}^{n}, B_{2}^{n}, \ldots$
and  that 
$$r_{i} - {1  \over n} < \mu(B_{i}^{n}) < r_{i}. $$
By part
(a) we may partition the open set $C:= B \backslash \bigcup_{i}B_{i}^{n}$
into clopen subsets $C_{1}, C_{2}, \ldots$ 
of measure less than $\ep$. If $\ep$ is 
sufficiently small then by taking finite unions of the sets $C_{i}$
we may find disjoint clopen subsets $D_{1}, D_{2}, \ldots$ of $C$,
all but
finitely many empty, such that 
 if we set
$B_{i}^{n+1}= B_{i}^{n} \cup D_{i}$ then we have
$$
r_{i} - {1 \over n+1} < \mu(B_{i}^{n+1}) < r_{i}.
$$
After we have defined $B_{i}^{n}$ for all $i$ and $n$ we see that
the open sets 
$B_{i} = \bigcup_{n=1}^{\infty} B_{i}^{n}$ are disjoint and
$\mu(B_{i}) = r_{i}$.

\endproof

\proclaim Definition. A { \bf column} $\T$ in $X$ consists of an ordered
collection of 
disjoint clopen sets of
the form $\set{B_{0},B_{1}, \ldots, B_{h-1}} $
together with integers
$$0= n_{0},n_{1}, \ldots, n_{h-1}$$ such that 
$T^{n_{i}}B_{0} = B_{i}$. 

\def\O{{\cal O}}

We denote this column by 
$\T=(B_{i},n_{i})_{i=0}^{h-1}$; $h$ is the {\bf height}
of $\T$ and $B_{0}$ is the { \bf base}. We will use the notation 
$\abs{\T} = \bigcup_{i} B_{i}$. The sets $B_{i}$ are the {\bf levels} of $\T$.
The {\bf width} of $\T$ is 
the measure of its base.
If $x \in \abs {\T}$ then
$\T x $ will denote  the level of  $\T$ which contains $x$. 
If $\T x = B_{i}$ the set 
$$\O_{\T}x = \set{T^{n_{j}}( T^{-n_{i}}x):j= 
0, \ldots ,h-1}$$
will be called the $\T$-{ \bf fiber} of $x$.
Roughly speaking what we call a column is called a tower in [HK] but 
we will reserve the term tower for the usual notion of a
Rohlin tower, that is, a  column
$\T=(B_{i},n_{i})_{i=0}^{h-1}$ with $n_{i} = i$.

If $T$ has no rational spectrum then the sets $B_{i}$ determine the 
integers $n_{i}$ but in general both need to be specified. 
Nonetheless we will sometimes refer to  $\{B_{0}, \ldots, 
B_{h-1}\}$ as a column with the understanding that this means
that the
${n_{i}}$ are clear from the
context, or, that we are asserting the existence of suitable $n_{i}$.

A {\bf slice} of  $\T=(B_{i},n_{i})_{i=0}^{h-1}$ is 
any column of the form $\T'=(E_{i},n_{i})_{i=0}^{h-1}$ such that 
$E_{i} \subset B_{i}$ . $\T'$ is 
determined by specifying the set $E_{0} \subset B_{0}$ and we will
call $\T'$ the slice of $\T$ over $E_{0}$.

\proclaim Definition. An {\bf array} $\T$ of height $h$ is a finite or 
countable 
collection of pairwise disjoint  columns $\T_{i}$ of the same height $h$.

\def\L{{\cal L}}

We will write $\abs \T = \bigcup_{i}\abs{\T_{i}}$.
Every column is also an array. The {\bf base} of $\T$ is the union of 
the bases of its columns and its width is the measure of its base.
The {\bf levels} of $\T$ are the levels of its columns and we denote the
set of levels by $\L(T)$. (Warning: the base of $\T$ is not a level of
$\T$ unless $\T$ has just one column.)
If the array $\T$ is contained in a measurable subset $E$ of $X$ we  say
it is an {\bf array partition }
of $E$ if its levels fills out $E$.
 Note that the total measure of $\T$ is
$hw$ where $h$ is the height of $\T$ and $w$ the width.
If $\T$ is an array  partition of $X$ we call 
it simply an {\bf array partition}.

\proclaim Definition. A {\bf sub-array} of an array $\T$ is an array $\cal 
T'$ such that each column of $\T'$ is a slice of  a column of $\T$.
If $\T'$ fills out $\T$ we call it a {\bf refinement} of
$\T$. 

\proclaim Lemma 3. Suppose $\T$ is an array of width $r$ and 
$r = \sum_{i=1}^{\infty} r_{i}$. Then $\T$ has disjoint subarrays $\T_{1}, 
\T_{2}, \ldots$ such that the width of $\T_{i} $ is $r_{i}$.

\proof Suppose the bases of the columns of $\T$ are
$B_{i}, i \in \N$. The base $B = \bigcup_{i}B_{i}$ is an
open set so there exist open subsets $C_{1},C_{2}, \ldots$ of $B$
such that $\mu(C_{i}) = r_{i}$. Each $C_{i}$ is a countable union of
disjoint clopen subsets $C_{i,j}, j \geq 1$. Let 
$\T_{i}$ be the subarray of $\T$
whose columns have bases $B_{k}\cap C_{i,j}, j,k  \geq 1$.

\endproof

Suppose $\T_{1} = \{B_{0}, \ldots , B_{h-1}\}$ and 
$\T_{2}=  \{C_{0}, \ldots , C_{k-1}\}$
are disjoint
columns of equal width. Suppose further that there is an integer 
$m$ such that $T^{m}B_{0} = C_{0}$. Then we can form the column 
$$\T = \{B_{0}, \ldots , B_{h-1},C_{0}, \ldots , C_{k-1}\}$$
called the 
{\bf concatenation} of $\T_{1}$ and $\T_{2}$. This notion can be 
extended in the obvious way to define the concatenation of columns
$\T_{1}, \ldots, \T_{s}$ provided their bases themselves form a column.
We will denote this concatenation by 
$\T=[\T_{1}, \ldots , \T_{s}]$.

\def\ht{\hat \T}

\proclaim Definition. An array $\ht$ is an {\bf extension} of an array
$\T$ if there is a refinement $\T'$ of $\T$ such that each column of 
$\ht$ is a concatenation of columns of $\T'$ and $\ht$ fills out
 $\T$. $\T'$ will be called the refinement of $\T$ 
associated to $\ht$. 

Note that $\ht$ and $\T'$ have the same set of levels. We observe that
a refinement is also an extension in a trivial way. If 
$\ht$ is an extension of $\T$ there is a natural projection from the 
set of levels of $\ht$ to those of $\T$ which we denote by $\pi_{\T}$:
$\pi_{\T} L = L'$ if $L \subset L'$.

\def\C{{\cal C}}

\proclaim Definition.
 Suppose $\T_{0},\T_{1},\ldots, \T_{m-1}$ are disjoint 
 arrays of equal widths.
 An array $\T$ will be called a {\bf stacking}
of   $\T_{0},\T_{1},\ldots, \T_{m-1}$ if there are refinements
$\T'_{0},\T'_{1},\ldots, \T'_{m-1}$
 such that each 
column of $\T$ is a concatenation $[\C_{0}, \ldots, \C_{m-1}]$, where
each
$\C_{i}$ is a column of $\T'_{i}$, and $\T$ fills out 
$\bigcup_{i}\abs{\T_{i}}$. Note that the order in which the 
$\T_{i}$ are listed is an essential part of this definition.

Evidently any extension of $\T$ is obtained as a stacking of 
$\T_{1}, \ldots , \T_{m}$ where  the $\T_{i}$ are disjoint subarrays of
$\T$ which fill out $\T$.
\proclaim Lemma 4. 
\item{(a)} If  $\T_{0},\T_{1},\ldots, \T_{m-1}$ are arrays of equal
width 
then there is a stacking $\T$ of  $\T_{0},\T_{1},\ldots, \T_{m-1}$.
\item{(b)}
If 
$\T$ is an array of height $h$ then $\T$ has an extension of height 
$hm$ for any $m \geq 1$.

\proof \item{(a)} Without loss of generality we assume $m=2$, say
$\T_{0} = {{\cal R}}$ and $\T_{1} = \S$, with bases $A$ and $B$.
By Lemma 1  there are 
disjoint 
clopen sets $A_1,A_2, \ldots \subset A$
 which fill out $A$ and 
integers $n_1, n_2, \ldots$ such that 
the sets $B_{i} := T^{n_{i}}A_{i}$
are disjoint and 
contained in $B$.
Let $\set{C_{i}}
$ and $\set{D_{i}}$ be the bases of 
the columns of ${{\cal R}}$ and $\S$ and set 
$$A_{ijk} = A_{i} \cap C_{j} \cap T^{-n_{i}}D_{k}$$
 and
$B_{ijk}= T^{n_{i}}A_{ijk}$. Then $\set{A_{ijk}}$ and $\set{B_{ijk}}$ 
are disjoint families of 
clopen subsets of $A$ and $B$ filling out
 $A$ and $B$. Let ${{\cal R}}_{ijk}$ and $\S_{ijk}$ be the slices of ${{\cal R}}$ and $\S$
 with bases
 $A_{ijk}$ and $B_{ijk}$. By construction the
 concatenations $[{{\cal R}}_{ijk},\S_{ijk}]$ exist. These concatenations are
 the
 columns of the desired stacking of ${{\cal R}}$ and $\S$.
    
\smallskip   

\item{(b)} If $\T$ has width $w$
simply divide $\T$ into subarrays $\T_{i}$, $i= 0, \ldots, m-1$ of 
width $w/m$ and then  use part (a) to stack these.

\endproof

\def\diam{\rm diam\ }

We define the {\bf diameter} of an array $\T$,
denoted $\diam \T$,  as the supremum of the
diameters of its levels.

\proclaim Lemma 5. Any array $\cal T$ has a refinement $\cal T'$
with  $\diam \T < \ep$.

\proof Without loss of generality
$\cal T$ has just one column $(B_{i},n_{i})_{i=0}^{h-1}$.
Using the fact that
$X$  is fractured it is easy to see that there is a countable partition $P$ 
of $X$ into clopen sets of diameter less than $\ep$. For 
$p = (p_{0}, \ldots , p_{h-1}) \in P^{h}$
let
 $$B_{p} = \bigcap_{i=0}^{h-1}T^{-n_{i}}(p_{i} \cap B_{i}).$$ Then the 
clopen
sets $B_{p},\  p \in P^{h}$  are disjoint and cover $B_{0}$.
The slices of $\T$ over the  sets $B_{p}, p \in P^{h}$
form the desired refinement.
\endproof

The following lemma is just the usual Rohlin lemma but  with a tower 
whose levels are clopen.

\proclaim Lemma 6. For every $h \in \N$ and $\ep > 0$ 
there is a clopen set $B$ such that the sets $B,TB, \ldots, T^{h-1}B$ 
are disjoint and cover at least $1-\ep$ of $X$.

\proof  Just repeat the usual proof of the Rohlin Lemma
(see for example [F]) starting from a 
small clopen set. Alternately one can take a measurable Rohlin tower,
approximate its base $B$ by a clopen set $C$ and then disjointify the
images 
of $C$.

\endproof

The following lemma is key to the proof of Theorem 2. It says
roughly that any
array partition has an extension whose fibers fill out large segments 
of the orbits of $T$.

\def\hT{\hat \T}

\proclaim Lemma 7. Suppose $\cal T$ is an array 
partition of $X$ 
and $\ep > 0$. Then one can find 
an extension $\ht$  of $\T$ and
a clopen set $X^{\sharp} $ such that $X^{\sharp} \subset 
\abs \hT$, $\mu(X^{\sharp}) > 1 - \ep$ and for all $x \in X^{\sharp}$ 
we have $Tx \in \O_{\hT}x$.

\def\L{{\cal L}}

\proof Given a column $\C= (C_{i},n_{i})_{i=0}^{h-1}$ we will 
refer to the integer $n 
= \max_{i,j}\abs{n_{i}-n_{j}}$  as the {\bf
spread} of 
$\T$.   Suppose $\T$ has height $s$ and columns $\T_{i}, i =1, 2,
\ldots$,
and let $\L = \L(T)$ 
denote the collection of levels of $\T$.
Fix a $\de >0$ to be specified  in the course of the proof and take
a finite collection $\set{\T_{i}, i = 1, \ldots, k}$ of columns of $\T$
which 
covers $1-\de$ of $X$ and agree to call these the {\bf good } columns.
Let $X' = \bigcup_{i=1}^{k}\abs{\T_{i}}$.
Let $m$ be the maximum of the spreads of the good  columns.

Find a clopen Rohlin tower for $T$ 
with base $E$ and height $h \gg m$,
to be further specified later, which covers more than $1-\de$ of 
$X$. For each $ { L} = (L_{0}, L_{1,},\ldots, L_{h-1}) \in \L^{h}$ let
$E_{L} = \bigcap_{j} T^{-j}(T^{j}E \cap L_{j})$. Then the sets
$E_{ L}, 
L\in \L^{h}$ are disjoint clopen subsets of $E$ which fill out $E$.
Each $E_{L}$ is the base 
of a  Rohlin tower of height $h$ and $T^{j}E_{ L} \subset L_{j}$
for $0 \leq j \leq  h-1$.

Let 
us call the union of the top and bottom $m$ levels of the tower the 
{\bf buffer}. We assume that $h$ has been chosen large enough
so that  the buffer has measure less than $\de$.
Fixing  for the moment any $L \in \L^{h}$ and 
any $j$ such that $ m < j <h-m$, let 
$ (C_{i},n_{i})_{i=0}^{s-1}$ 
(temporarily) denote the column of $\T$ which contains $T^{j}E_{ L}$.
If 
$T^{j}E_{L} \subset C_{r}$ 
then we will denote the slice of $ (C_{i},n_{i})_{i=0}^{s-1}$
over $T^{-n_{r}}T^{j}E_{L}$ by
$\C_{L,j}$, so that $T^{j}E_{ L}$ is the $r$-th level of $\C_{L,j}$.
We will be interested in those pairs $(L,j)$, call them {\bf good},
such that $m < j < h-m$
and
$T^{j}E_{L}$ is contained in a level of one of the good
columns  $\T_{i}, 1 \leq i \leq k$.
Note that  for a good $(L,j)$ 
the levels of $\C_{L,j}$ are all levels of the $T$-tower
of height $h$ over $E_{L}$, since they are of the form 
$T^{s}T^{j}E_{L}$ with $\abs s \leq m$ and $m < j < h-m$. Moreover,
for a fixed $L \in \L^{h}$ and any two good $(L,j)$ and $(L,j')$,
$\C_{L,j}$ and $\C_{{ L},j'}$
are either disjoint or equal. This is clear:
$\C_{L,j}$ and $\C_{L,j'}$ are both slices of $\T$ and each consists
of levels of the $T$-tower of height $h$ over $E_{L}$ so if they 
are not disjoint they share a level. Clearly if two slices of
$\T$ share a level they are identical.

  \def\G{{\cal G}}

  
Let $\S_{L}$ denote the set of distinct slices $\C_{L,j}, m < j < h-m$ 
such that $(L,j)$ is good.
Since $\mu(X') > 1 - \de$ and the measure of the buffer is 
less than $\de$, a `Fubini' argument
shows that $E$ is  $(1-\de_{1})$-covered by $E_{ L}$'s for which the
slices in
$S_{ L}$ cover at least a fraction
$1-\de_{1}$ of $\bigcup_{j}T^{j}E_{ L}$,
where $\de_1 = \sqrt{2\de}$.
Call these the good $L$'s and denote the set of good 
$L$'s by $\G$. Let $r$ be the least integer such that
$rs \geq (1- \de_{1}) h$. Then for each good $L$ the cardinality of
$\S_{L}$ is  at least $r$ so we  may choose  $r$ elements of 
$\S_{L}$ and concatenate them in any order to form a column $\C_{L}$ 
of height $rs$.
The concatenation is possible because the bases of these $r$ columns
themselves form a column,  as they are all levels of the
tower of height $h$ over $E_{L}$.
 Our choice of $r$ implies that if  $L$ is good then
 $\C_{L}$ covers at least a fraction 
 $1-\de_{1}$ of $\bigcup_{j=0}^{h-1}T^{j}E_{ L}$. Thus we have
 $$
\mu(\bigcup\set{\abs{\C_{L}}: L \in \G}) > 
(1-\de_{1})^{2}(1-\de) =: 1-\de_{2}.
 $$

Choose a finite subset $\G' \subset \G$ such that the measure of the clopen
set
$$
X^{\flat} =\cup\set{\abs{\C_{L}}: L \in \G'})
$$
 is greater then $1-2 \de_{2} =: 1 -\de_{3}.$
 The $\C_{L}, L \in \G'$, form an array $\T_{0}$ contained in 
 $X'$  with finitely many columns which will be
 part of the desired extension $\ht$. Since $\abs{\T_0}$ is clopen
 what remains of $\T$ is again an array $\T_{1}$. It remains only
 to form any extension $\T_{2}$ of $\T_{1}$ of height $rs$ and adjoin it
 to $\T_{0}$ to obtain  $\ht$.

 Now set 
 $X^{\sharp} = X^{\flat} \cap T^{-1}X^{\flat} \cap (T^{h-1}E)^{c}$
 so $$\mu(X^{\sharp}) > 1-2\de_{3} - {1 \over h} > 1 - \ep,$$
if $\de$ is sufficiently small and $h$ sufficiently large.
 If $x \in X^{\sharp}$  then $x$ and $Tx$ belong to
 the same fiber of the $T$-tower of height $h$ over $E$, since
 $x \notin T^{h-1}E$. Since both $x$ and $Tx$ belong to $X^{\flat}$ 
 they  are in the same fiber of $\C_{L}$ for some $L \in  \G$, hence
 $x$ and $Tx$ are in the same fiber of $\hT$.

 \endproof

\def\hs{\hat{\S}}
\def\ht{\hat {\T}}

\proclaim Definition.
Suppose that $\cal T$ and $\cal S$ are array partitions of $X$ and
$Y$ of 
the same height. An {\bf array 
map} from $\cal T$ to $\cal S$ is a map from the
set of levels of $\T$ to those of $\S$
 which maps individual columns bijectively 
 to individual columns, preserves the order on each column and maps 
 $\mu$ to $\nu$: if $L \in \L(S)$ is a level of  $\cal S$ 
then $\mu ( \bigcup \ph^{-1} \set{L}) = \nu (L)$.

To clarify the notation, the expression  $ \bigcup \ph^{-1} \set{L}$ can be read literally but
it could also be written as $\bigcup\set{L' \in \L(\T): \ph(L') = L}$.
The map 
$\ph$ can also be viewed as a mapping from the set of  columns of $\T$ to 
the set of columns of $\S$.
One may visualize $\ph$ as 
lumping together a sub-array of  $\T$
(whose columns are columns of $\T$ rather than  slices)
with total width $w$ to give one column
of $\S$ of width $w$. We want to stress however that $\ph$ is a 
mapping of sets, not points. 

Suppose $\ph : \cal T \to \cal S$ is an array map 
and  suppose $\ht$ and $\hs$ are extensions of $\T$ and $\S$ 
with associated refinements $\T'$ and $\S'$. We will use $\pi_{\T}$
to 
denote the projection from $\ht$ to $\T$ and $\pi_{\S}$ for the
projection from 
$\hs$ to $\S$.
An array map $\ps: 
\hs \to \ht$ will be called an {\bf extension} of  $\ph$ if
$\ph \circ \pi_{\T} \circ \ps = \pi_{\S}$. The fact that 
$\ps$ preserves the order on columns implies that
$\ps$ is also an array map from $\S'$ to $\T'$.
Similarly if $\psi: \ht \to \hs$ we call it
an extension of $\ph$ if $\pi_{\S} \circ \psi = \ph \circ \pi_{\T}$ 
and again $\ps$ is automatically also an array map from $\T'$ to
$\S'$. Analogous definitions hold if $\ph : \S \to \T$. 
By chasing commutative  diagrams it is easy to see  that if
$\T''$ and $\S''$ are extensions of $\T'$ and $\S'$ and
$\chi: \T'' \to \S''$  is an extension of $\psi$ then it is an extension
of $\ph$, whatever the directions of the maps $\ph$, $\psi$ and $\chi$.

\proclaim Lemma 8 (Copying Extensions). Suppose $\ph : \cal T \to \cal S$ is 
an array map 
and  suppose $\ht$ 
is an extension of $\cal T$. Then there is 
an extension $\hs$ of $S$ and an array map $\psi:\hs \to \ht$ 
which extends $\ph$.

\proof  Suppose the heights of $\T$ and $\ht$ are $h$  and $mh$. 
Let $\T'$ denote the refinement of $\T$ associated to $\ht$.
Let $\set{\S_{i}:i \in I}$ denote the set of columns of $\S$. Label
each column
of $\T$ by  $i \in I$ if it is mapped to $\S_{i}$.
For each column $\C \in  \ht$ let $w_{\C}$  denote its width and
${ i}_{\C} = ({ i}_{\C,0}, \ldots ,i_{\C,m-1})$  its string of labels:
$\C$ is a concatenation $[\C_{0}, \ldots , \C_{m-1}]$ where each 
$\C_{j}$ is a slice of a column of $\T$ which has label $i_{\C,j} $.

Now find a system of disjoint subarrays 
$$\set{\S_{\C,j}, \C \in  \ht, 0\leq j < m}$$
of $\S$
such that $\S_{\C,j}$ is a subarray of
$\S_{i_{\C,j}}$ and has width $w_{\C}$.  This is possible because of
Lemma 3 and the
measure-preserving character of $\ph$.  Observe that the $\S_{\C,j}$
 fill out  $\S$.  Now use Lemma 4 to
find a stacking $\S_{\C}$ of
$\S_{\C,0},\S_{\C,2}, \ldots
\S_{\C,m-1} $.  The
$\S_{\C},\ \C \in \ht$,
 together form the array partition $\hs$ we seek.  Now define
$\psi$ by mapping each column of $S_{\C}$ to $\C$.
         
\endproof

\bigskip
We are now ready to prove Theorem 2.
 Using  Lemmas 5,7 and 8  we inductively 
construct a sequence of arrays $\T_{i}, i = 0,1, \ldots$ in $X$, a 
 sequence 
of arrays $\S_{i}, i =0,1, \ldots$ in
 $Y$ and array maps $\ph_{i}: \T_{i} \to 
\S_{i}$ for even $i$ and $\ps_{i}: \S_{i} \to \T_{i}$ for odd $i$ 
such that, for each $i$,  $\T_{i+1}$ extends $\T_{i}$ and  $\S_{i+1}$ extends 
$\S_{i}$  and
\smallskip

\item{(i)} for each even $i$, $\ph_{i}$ extends $\psi_{i-1}$ and
$\psi_{i+1}$ extends $\ph_{i}$.
\smallskip
\item{(ii)}  $\diam\T_{i} \to 0$ and $\diam \S_{i} \to 0$
\smallskip
\item{(iii)} for each even $i>0$ 
there is a clopen set $X_{i} \subset \T_{i}$ such that
$\mu(X_{i}) > 1 - {1 \over i}$ and 
for  each $x \in X_{i}$ we have $Tx \in 
\T_{i}x$;   for each even $i$ there is a clopen
$Y_{i} \subset \S_{i}$ satisfying the analogous conditions.

\smallskip
We remark that (i) implies that any $\ph_{i}$ or $\psi_{i}$ is an
extension of any $\ph_{j}$ or $\ps_{j}$ for $j < i$.
Let $X' = \bigcap_i \abs{\T_{i}}$ and $Y' = \bigcap_i \abs{\S_{i}}$. These are
$G_{\de}$ subsets of full measure, not necessarily invariant.
We define a map $\ph: X' \to Y'$ as follows. For $x \in X'$, $\set{\T_{2i}x}$
is a decreasing sequence of clopen sets whose intersection is 
$\set x$.
Since $\ph_{2i+2}$ extends $\ph_{2i}$ it follows that 
$\ph_{2i}(\T_{2i}x)$ is also a decreasing sequence of clopen sets.
Since the diameters go to $0$ it follows that the intersection
is a singleton
$\set y$ such that $y \in Y'$ and we define $\ph(x)= y$.
We define $\ps: Y' \to X'$ similarly using the maps
$\ps_{2i+1}$. Note that the definition of $\ph$ boils down to
$\S_{2i}\ph(x) = \ph_{2i}(\T_{2i}x)$.

\def\diam{{{\rm diam} }}
\def\id{{\rm id}}

\def\L{{\cal L}}
Let ${\cal L} = \bigcup_i \L(\S_{2i})$. The collection
$ \set{L \cap Y': L \in \L}$ generates the topology of
$Y'$, since its members 
are open and $ \diam \S_{2i} \to 0$.
Thus to show $\ph$ is
continuous it suffices to show that for all $i$ and all  $L \in
\S_{2i}$,
we have $\ph^{-1}(L \cap Y')$ open in the relative 
topology of $X'$.  But it is clear that
$$\ph^{-1} (L \cap Y') = (\bigcup \ph^{-1}_{2i} \set{L}) \cap X'\ :$$
evidently 
$$\ph  (( \bigcup \ph^{-1}_{i} \set{L}) \cap X') \subset L \cap Y'$$
and on the other hand if
$x \in X'$  and $\T_{2i}x$ is disjoint from  $\bigcup 
\ph^{-1}_{2i} \set{L} \cap X'$ then $\ph x$ is in
$\ph_{2i}(\T_{2i}x)$ which is disjoint  from $L$. Thus
$\ph^{-1}(L \cap Y')$ is open in $X'$,
we  have shown that $\ph$  is continuous, and the same
argument holds for $\ps$.

Next we check that $\ph $ and $\psi$ are inverse maps. 
 For $x \in X'$ we have
$$
\eqalign{
\T_{2i-1}\ps\ph x = & \psi_{2i-1}\S_{2i-1}\ph x \cr= &
\psi_{2i-1}\pi_{\S_{2i-1}}\S_{2i}\ph x \cr  =&
\ps_{2i-1}\pi_{\S_{2i-1}} 
\ph_{2i}\T_{2i}x \cr = &
\pi_{\T_{2i-1}}\T_{2i}x \cr = &\T_{2i-1}x.
}
$$
Since  $\diam \T_{2i-1}x \to 0$  it follows that 
$\ps \ph x = x$. 
 Similarly $\ph \circ \ps = \id_{Y'}$ so $\psi$ is a homeomorphism
 from $X'$ onto $Y'$ and $\psi = \ph^{-1}$.

Now we claim that  for $x \in X'$ we have
$\ph (\O_{\T_{i}}x) = \O_{\S_{i}}\ph(x)$ for 
 each even $i$. (Note that if $x \in X'$ then $\O_{\T_{i}}x \subset X'$
 for each $i$).
 For any even $j > i$ let
 $\T_{ij}$ denote the refinement of $\T_{i}$ determined by its 
extension $\T_{j}$ and $\C_{ij}x$ the column of $\T_{ij}$ containing 
$x$, with analogous definitions in $Y$. 
Clearly $\bigcap_{j}\abs{\C_{ij}x} = \O_{\T_{i}}x$.
Since $\ph_{j}$
is an array map from $\T_{ij}$ to $\S_{ij}$ it is clear that
$\ph(\abs{\C_{ij}x} \cap X') \subset \abs{\C_{ij}\ph(x)}$. 
Thus 
$$
 \ph(\O_{\T_{i}}x)= \ph(\bigcap_{j}\abs{\C_{ij}x}) =
 \bigcap_{j}\ph(\abs{\C_{ij}x })
\subset \bigcap_{j}\abs{\C_{ij}\ph(x)} = \O_{\S_{i}}(\ph x).
$$
Since $\ph$ is injective we conclude that $\ph(\O_{\T_{i}}x) = 
\O_{\S_{i}}x.$

Recall that $X_{n}$, as defined in (iii) above, is clopen.
Let $X^{*} =  (\bigcup _{i}X_{2i}) \cap X' \cap T^{-1}X'$, so $X^{*}$ is 
a $G_{\de}$ of full measure.
We claim 
that for
$x \in X^{*}$ we have $\ph(Tx) \in \O_{S}\ph(x)$.
This is clear:
if $x \in X^{*} $, so $x \in X_{2i}$ for some $i$,
then $Tx \in \O_{\T_{2i}}x$ so 
$$
\ph (Tx) \in \ph (\O_{\T_{2i}}x) = 
\O_{\S_{2i}} \ph x \subset \O_{S} \ph x
$$
as we claimed.
Similarly  we obtain a dense $G_{\de}$ of full measure $Y^{*} 
\subset Y' \cap S^{-1}Y'$ such that for $y\in Y^{*}$ we have 
$\psi(Sy) \in \O_{T}\psi y $.

Now  $X^{*}$ and $Y^{*}$ need not be invariant nor do we have
$\ph(X^{*}) = Y^{*}$.
We remedy this as follows. 
  Given a $G_{\de}$ subset $E$ of full measure in $X$ we will write
$ E^{T}$ for the set $\bigcap_{n \in \Z}T^{n}E$, which is
an invariant $G_{\de}$ 
subset of $E$ of
full measure.
We will use the same notation in $Y$. 
Let $$X_{1}= (X^* \cap \psi Y^{*})^{T},$$
$$Y_{1} =  (\ph(X_{1}))^{S},$$
$$X_{2}=( \psi(Y_{1}))^{T},
$$ $$Y_{2} =  (\ph (X_{2}))^{S}$$ 
and so on. Let $X_{\infty} = \bigcap_{n}X_{n}$ and $Y_{\infty} = 
\bigcap_{n} Y_{n}$. Then $X_{\infty}$ and $Y_{\infty}$ are invariant 
$G_{\de}$ subsets 
of $X^{*}$ and $Y^*$, both having full measure, and evidently 
 $\ph X_{\infty} = Y_{\infty}$. For $x \in X_{\infty}$
 $\ph(Tx) \in \O_{S} \ph(x)$  since $ X_{\infty} \subset X^{*}$ and 
 because $X_{\infty}$ is invariant
it follows 
easily that
$\ph \O^{+}_{T}x \subset \O_{S}\ph(x)$ for all $ x \in X_{\infty}$,
where $\O^{+}$ denotes the forward orbit.
Now since we also have 
$$\ph(\O^{+}_{T}x) \subset \ph(\O^{+}_{T}
T^{-1}x) \subset 
\O_{S}\ph(T^{-1}x)$$
 it follows that 
$\O_{S}\ph(x) =\O_{S}\ph(T^{-1}x) $
and in particular $\ph(T^{-1}x) \in \O_{S}\ph(x)$. In a similar way 
we find that $\ph(T^{-i}x) \in \O_{S}\ph(x)$ for all $i > 0$ so we 
conclude that $\ph(\O_{T}x) \subset \O_{S}\ph(x)$ for any $x \in X_{\infty}$.
All we used to show this is that $X_{\infty}$ is invariant and 
contained in  $X^{*}$. It follows that the corresponding fact for $\psi$
holds as well, even though the definitions of $X_{\infty}$ and 
$Y_{\infty}$ are not symmetric.
Using this we see that
$$
\O_{T}x = \psi  \ph \O_{T}x \subset \psi \O_{S} \ph(x) 
\subset \O_{T} \psi(\ph(x)) = \O_{T} x.
$$
It follows that all the containments are in fact equalities and in 
particular that $\ph (\O_{T}x) = \O_{S} \phi(x)$. This concludes the
proof of Theorem 2.

\bigskip

\noindent{Section 3.  The II$_\infty$ case.}
\smallskip
In this section we shall prove the following result by reducing it
to Theorem 1.

\def\gd{$G_\de$\ }

\proclaim Theorem 3. 
Suppose $X$ and $Y$ are Polish spaces with 
non-atomic infinite $\si$-finite measures
 $\mu$ and $\nu$. Suppose that $X$ and $Y$ 
have at least one open subset
of finite non-zero measure.
 Suppose that $T$ and $S$ are ergodic 
measure-preserving homeomorphisms of $(X,\mu)$ and $(Y,\nu)$.
Then there are invariant $G_{\de}$ subsets $X' \subset X$ and $Y' 
\subset Y$ of full measure and a homeomorphism $\ph: X' \to Y'$ which 
maps $\mu|_{X'}$ to $\nu|_{Y'}$ and maps $T$-orbits onto $S$-orbits.

\proof
As in the measurable case,
we shall prove Theorem 3 by inducing on a set of
finite measure. We begin with some remarks about inducing in our setting.
First observe that the hypotheses of Theorem 3 imply that $T$ is
conservative and that for any non-null set $E \subset X$ the set
$\bigcup_{i>0}T^{i}E$ is co-null.
If $U$ is a non-empty clopen subset of $X$ then it is easy to see that
the  induced transformation $T_U$ is a homeomorphism from an open subset
$U'$ 
of $U$, such that $\mu(U \backslash U') = 0$,
to another open subset of $U$. If we assume further that for 
all $x \in X$ the
$T$-orbit  of $x$ intersects $U$ infinitely often in both positive and
negative time then $T$ will be a homeomorphism from $U$ to itself.
Call such sets $U$ {\bf $T$-good}. 
If $U$ is not $T$-good, observe that the set $X'$ of points in 
$X$ whose
$T$-orbit 
intersects $U$ infinitely often in both positive and
negative time
is an invariant co-null $G_{\de}.$ Setting
$U' = U \cap X'$ then $U'$ is clopen in $X'$ and it is
$T_{X'}$-good. (The induced map $T_{X'}$ is the same as
the restriction of $T$ to $X'$). This means that for the
purpose of proving Theorem 3 there is
no loss of generality in assuming that $U$ itself is $T$-good. Note also
that if $U$ is $T$-good and $X'$ is any invariant co-null \gd then
$U' = U \cap X'$ is $ T_{X'}$-good. Finally, any clopen superset
of a $T$-good set is $T$-good.

To prove Theorem 3 we may assume,
 as in the proof of Theorem 1, that $X$ and $Y$ are fractured and
that $\mu$ and $\nu$ have full support. Let $U$ be a
\nec of finite measure in $X$.
Then $\bigcup_{i \geq 0}T^{i}U  $ is co-null \gd
and after restricting to  $(\bigcup_{i \geq 0}T^{i}U)^{T}  $
we might as well assume that  $\bigcup_{i \geq 0}T^{i}U  = X $.
 It follows easily that $X$ can also be expressed as a countable 
disjoint union of
 of  clopen sets $U_i$ of finite measure.
By the above remarks there is no loss of generality in assuming
that  each $U_i$ is $T$-good. Applying Lemma 2(b) to the induced system
$(U_i, \mu|_{U_i}, T_{U_i})$ we see that $U_i$ contains disjoint
open subsets $U_{i,j}$ of any desired measures which sum to $\mu(U_i)$.
From this it follows that we can find disjoint open subsets
$E_0, E_1, \ldots$ of $X$ such that $\mu(E_i) = 1$ for all $i$ and
$X':=\bigcup _{ i \geq 0} E_i$ is co-null. Restricting to $X'$ there is
no loss of generality in assuming that $\bigcup _{ i \geq 0} E_i = X$.
In addition we may assume that  $E_0$ is $T$-good. Since
$E_0 \cup E_i$ is also $T$-good we may apply Lemma 1 to $T_{E_0 \cup E_i}$
to find open subsets $E_0'$ and $E_i'$ and a homeomorphism $\th_i$ from
$E_0'$ onto $E_i'$ such that $\th_i$ is a piecewise power of 
$T_{E_0 \cup E_i}$, and hence also a piecewise power  of $T$.
Cutting down to $( \bigcup E_i)^{T}$ we may assume that $\th_i$ maps all
of $E_0$ to all of $E_i$. It is easy to see that
$\th_{i} (\ot x \cap E_{0}) = \ot x \cap E_{i}$ for any
$x \in X$. Letting $T_{0} = T_{E_{0}}$ we observe also that 
for $x \in E_{0}$ we have ${\cal O}_{T_{0}} x = \ot x \cap E_{0}$.

Now go through the same process in $Y$ to obtain corresponding sets
$F_i$ and homeomorphisms $\psi_i$. Let $T_0 = T_{E_0}$,
$S_0 = S_{F_{0}}$, $\mu_{0} = \mu|_{E_{0}}$ and $\nu_{0}= \nu|_{F_{0}}$.
Applying Theorem 1 to the homeomorphisms $T_{0}$ and $S_{0}$ with 
measures $\mu_{0}$ and $\nu_{0}$ we find invariant
\gd subsets $E_{0}' \subset E_{0}$ and $F_{0}' \subset F_{0}$, both 
of measure one and a homeomorphism $\ph : E_{0}' \to F_{0}'$ which 
maps $T_{0}$-orbits onto $S_{0}$-orbits and $\mu_{0}$ to
$\nu_{0}$. Now set
 $E_{i}'= \th_{i}E_{0}'$, $X' = \bigcup E_{i}'$,  $F_{i}' = \psi_{i}F_{0}'$ and
 $Y' = \bigcup F_{i}'$. Then $X'$ and $Y'$ are invariant, co-null \gd 
 subsets of $X$ and $Y$. Define a homeomorphism $\Phi:X' \to Y'$ by
 $$
 \Phi(x) = \psi_{i}\ph \th_{i}^{-1}x \hbox{ for } x \in E_{i}'.
 $$
 Then it is clear that $\Phi$ is the desired orbit equivalence.
 \endproof

\noindent{References}
\bigskip

\noindent [D1]{} H. Dye, {On groups of measure-preserving
transformations I}, Amer. J. Math., 81 (1959), 119-159.
\smallskip
\noindent [D2] H.Dye, {On groups of measure-preserving
transformations II},  Amer. J. Math.  85  (1963) 551-576.
\smallskip
\noindent [F] N. Friedman, { Introduction to ergodic theory.},
Van Nostrand Reinhold Mathematical Studies, No. 29. (1970).
\smallskip
\noindent [GPS] T. Giordano, I. Putnam and C. Skau, 
Topological orbit equivalence and $C\sp *$-crossed products. 
J. Reine Angew. Math.  469  (1995), 51-111.
\smallskip
\noindent [HIK] 
A. Hajian, Y. Ito and S. Kakutani,
{Full groups and a theorem of Dye},
Advances in Math. 17 (1975), 48-59. 
\smallskip
\noindent  [HK]{} T. Hamachi and M. Keane, {Finitary orbit equivalence of
odometers}, Bull. London Math. Soc., 38(2006), 450-458.
\smallskip
\noindent  [HKR]{} T. Hamachi, M.S.Keane and M.K. Roychowdhury, 
 {Finitary orbit equivalence and measured Bratteli diagrams} (preprint).
\smallskip
\noindent [HO] T. Hamachi and M. Osikawa, Ergodic groups 
of automorphisms and Krieger's theorems, Seminar on Mathematical
 Sciences, 3, Keio University, Department of Mathematics, Yokohama, 1981.
\smallskip
\noindent  [J] 
A. del Junco,
Bernoulli shifts of the same entropy are finitarily and unilaterally 
isomorphic,
Ergodic Theory Dynam. Systems 10 (1990),  687-715. 
 \smallskip
\noindent [KW] 
Y. Katznelson and B. Weiss,
The classification of nonsingular actions, revisited,
Ergodic Theory Dynam. Systems 11 (1991),  333-348. 
\smallskip
\noindent  [KS1] M. Keane and M. Smorodinsky, Bernoulli schemes of the same
 entropy are finitarily isomorphic,
 Ann. of Math. (2)  109  (1979), 397-406.
 \smallskip
 
\noindent [KS2] M. Keane and M. Smorodinsky,
{Finitary isomorphism of irreducible Markov shifts},
Israel J.Math., 34 (1979), 281-286.
\smallskip
\noindent [K] W. Krieger,  On ergodic flows and the isomorphism of factors,
  Math. Ann.  223  (1976),  19-70.
\smallskip
\noindent  [R1]{} M.K. Roychowdhury, {$\{m_n\}$-odometer and the binary
 odometer are finitarily orbit equivalent} 
 (To appear, Contemporary Mathematics, AMS,
 Edited by I. Assani, UNC Chapel Hill).
\smallskip
\noindent  [R2]{} M.K. Roychowdhury, {Irrational rotation of the circle and
 the binary odometer are finitarily orbit equivalent},
 Publ. RIMS, Kyoto Univ., 43 (2007), 385-402.
 
\smallskip
\noindent [RR1]{} M.K. Roychowdhury and D.J. Rudolph, 
{All uniform odometers are finitarily Kakutani equivalent} 
(preprint).
\smallskip
\noindent  [RR2]{} M.K. Roychowdhury and D.J. Rudolph,
 {Any two irreducible Markov chains of equal entropy
 are finitarily Kakutani equivalent}
(To appear, Isr. J. of Math).
\smallskip
\noindent  [RR3]{} M.K. Roychowdhury and D.J. Rudolph, 
 {Any two irreducible Markov chains  are finitarily orbit equivalent}
 (preprint).
\smallskip
\noindent [SSW] 
D. Sullivan, B. Weiss and J. D. Wright,
Generic dynamics and monotone complete $C\sp *$-algebras,
Trans. Amer. Math. Soc. 295 (1986),  795-809. 
\smallskip

\end